\documentclass[11pt]{article}
\usepackage{amsmath,amssymb,amsfonts}

\newtheorem{theorem}{Theorem}
\newtheorem{proposition}{Proposition}
\newtheorem{corollary}{Corollary}

\def\D{{\cal D}}
\def\R{{\mathbb R}}
\def\C{{\mathbb C}}
\def\Z{{\mathbb Z}}

\def\Re{{\mathrm{Re}\, }}
\def\Im{{\mathrm{Im}\, }}

\begin{document}

\title{Surfaces in the four-space and the Davey--Stewartson equations}
\author{Iskander A. TAIMANOV
\thanks{Institute of Mathematics, 630090 Novosibirsk, Russia;
e-mail: taimanov@math.nsc.ru}}
\date{}
\maketitle

\section{Introduction}
\label{sec1}

The Weierstrass representation for surfaces in $\R^3$ \cite{K1,T1}
was generalized for surfaces in $\R^4$ in \cite{PP} (see also \cite{FLPP}).
This paper  uses the quaternion language and
the explicit formulas for such a representation were written
by Konopelchenko in \cite{K2} for constructing surfaces which admit
soliton deformation governed by the Davey--Stewartson equations.
This generalizes his results from \cite{K1} where
he introduced the formulas for inducing surfaces in
the three-space which involve
a Dirac type equations and defined for such surfaces
a deformation governed by the modified Novikov--Veselov (mNV) equations.

It was shown in \cite{T1}
that the formulas for inducing surfaces in $\R^3$ \cite{K1}
describe all surfaces and that the modified Novikov--Veselov equation
deforms tori into tori preserving the Willmore functional
which naturally arises and plays an important role in this representation.
The spectral curve of the corresponding Dirac operator is
invariant under this deformation.
The global Weierstrass representation at least for real analytic
surfaces could be obtained by an analytic continuation from a local
representation. Thus a moduli space of immersed tori is
embedded into the phase space of an integrable system
with the Willmore functional and, moreover, the spectral curve
as conservation quantities.

Looking forward to understand the spectral curves for tori in $\R^4$ we
consider in this paper the analogous problems for surfaces in $\R^4$ and
show that this case is very different from the three-dimensional case,
in particular,
by the following features which were overlooked until recently:

\begin{itemize}
\item
for tori in $\R^4$ every equation from the Davey--Stewartson (DS)
hierarchy describes not one but
infinitely many geometrically different soliton deformations;

\item
the multipliers on
the spectral curve for a torus in $\R^4$ are not uniquely defined
and different complex curves in $\C^2$ (the spectral curves immersed via the
multipliers) are invariants of different
DS deformations.
\end{itemize}

The reason for that is quite clear and consists basically in the
nonuniqueness of a Weierstrass representation.

A surface in $\R^3$ is constructed in terms of one vector function $\psi$
(spinor) which is a lift of the Gauss mapping into non-vanishing spinors.
Such a lift is defined up to a sign by fixing a conformal parameter on
the surface. This function $\psi$ satisfies a Dirac equation.

A surface in $\R^4$ is constructed in terms of
two vector functions $\psi$ and
$\varphi$ which form again a lift of the Gauss mapping. However in this case
by fixing a conformal parameter one defines a lift only up
to a gauge transformation given by $e^f$ where $f$ is any smooth function.
Moreover not every lift satisfies the Dirac equations and the lifts meeting
these equations are defined up to gauge transformations $e^h$ where $h$ is a
any holomorphic function.

In particular, given a Weierstrass representation of  a surface
$\Sigma \subset\R^4$ and a domain $W \subset \Sigma$
we can replace a representation of the domain by
gauge-equivalent using  a transformation $e^h$ where $h$ is a
holomorphic function on $W$ which is not analytically continued onto the
surface. Thus we obtain a representation of a domain which is not continued
(i.e., expanded to a representation of a surface).
This also makes a difference with the case of surfaces in $\R^3$.

Another important point is that
the DS equations contain the additional potentials which are defined
by resolving the constraint equations. Such resolutions are
not unique and we
have to choose the potentials carefully to make the DS deformations
geometric:
for some special choices of the additional potentials the corresponding
DS deformations map tori into tori preserving the Willmore functional.
However in general this is not the case and we show how to achieve that in
\S \ref{sec4}.

The work was supported by RFBR (grant 03-01-00403) and
Max-Planck-Institute on Mathematics in Bonn.

We thank U. Abresch for comments and discussions.

\section{Explicit formulas for a representation and so\-liton deformations}
\label{sec2}

Let us recall the explicit formulas for inducing a surface
and its soliton deformation via the Davey--Stewartson equation.

The following proposition is derived by straightforward computations.

\begin{proposition}[\cite{K2}]
\label{first}
Let vector functions $\psi$ and $\varphi$ be defined in a
simply-con\-nec\-ted domain $W \subset \C$
(with a complex parameter $z$) and meet the Dirac equations
$$
\D \psi = 0, \ \ \
\D^\vee \varphi = 0
$$
where
$$
\D =
\left(
\begin{array}{cc}
0 & \partial \\
-\bar{\partial} & 0
\end{array}
\right)
+
\left(
\begin{array}{cc}
U & 0 \\
0 & \bar{U}
\end{array}
\right), \ \ \ \
\D^\vee =
\left(
\begin{array}{cc}
0 & \partial \\
-\bar{\partial} & 0
\end{array}
\right)
+
\left(
\begin{array}{cc}
\bar{U} & 0 \\
0 & U
\end{array}
\right).
$$
Then the $1$-forms
$$
\eta_k = f_k dz + \bar{f}_k d\bar{z}, \ \ \ k=1,2,3,4,
$$
with
\begin{equation}
\label{wr1}
\begin{split}
f_1 = \frac{i}{2} (\bar{\varphi}_2\bar{\psi}_2 + \varphi_1 \psi_1),
\ \ \ \
f_2 = \frac{1}{2} (\bar{\varphi}_2\bar{\psi}_2 - \varphi_1 \psi_1),
\\
f_3 = \frac{1}{2} (\bar{\varphi}_2 \psi_1 + \varphi_1 \bar{\psi}_2),
\ \ \ \
f_4 = \frac{i}{2} (\bar{\varphi}_2 \psi_1 - \varphi_1 \bar{\psi}_2)
\end{split}
\end{equation}
are closed and the formulas
\begin{equation}
\label{wr2}
x^k = x^k(0) + \int \eta_k, \ \ \ k=1,2,3,4,
\end{equation}
define a surface in $\R^4$ (here the integral is taken over any path
in $W$ and by the Stokes theorem does not depend on a choice of path).

The induced metric equals
\begin{equation}
\label{metric}
e^{2\alpha} dzd\bar{z} =
(|\psi_1|^2+|\psi_2|^2)(|\varphi_1|^2+|\varphi_2|^2)dz d\bar{z}
\end{equation}
and the norm of the mean curvature vector
${\bf H} = \frac{2 x_{z\bar{z}}}{e^{2\alpha}}$
meets the equality
\begin{equation}
\label{potential}
|U| = \frac{|{\bf H}| e^\alpha}{2}.
\end{equation}
\end{proposition}

For $U=\bar{U}$ and $\psi = \pm \varphi$ these formulas reduce to
the Weierstrass representation for surfaces in $\R^3$.

The existence of a local representation of any surface in $\R^4$
by these formulas is not proved in \cite{K2} although it was indicated in
\cite{PP} that the Weierstrass representation for surfaces in $\R^3$ is
generalized for surfaces in $\R^4$ and involves in this case two
vector functions $\psi$ and $\varphi$ and a complex valued potential $U$.

We expose such a derivation in the next section revealing some
features not taking place in the three-dimensional case.
Remark that for Lagrangean surfaces in $\R^4$ this representation was
discovered in other terms by Helein and Romon \cite{HR}.

Let
$$
L =
\left(
\begin{array}{cc}
0 & \partial \\
-\bar{\partial} & 0
\end{array}
\right)
+
\left(
\begin{array}{cc}
-p & 0 \\
0 & q
\end{array}
\right).
$$
Let us consider deformations of this operator which take the form of
Manakov's $L,A,B$-triple:
\begin{equation}
\label{triple}
L_t + [L,A_n] - B_n L=0
\end{equation}
or
$$
[L, \partial_t - A_n] + B_n L=0.
$$

Notice that if $L$ meets (\ref{triple}), then the solution of the equation
$$
L \psi = 0
$$
is evolved as follows:
$$
\psi_t = A_n \psi.
$$

The following two propositions are proved by straightforward computations.

\begin{proposition}
\label{second}
For
$$
A_2 =
\left(
\begin{array}{cc}
-\partial^2 - v_1 & q\bar{\partial} - q_{\bar{z}} \\
-p\partial + p_z & \bar{\partial}^2 + v_2
\end{array}
\right),
$$
$$
B_2 =
\left(
\begin{array}{cc}
\partial^2 + \bar{\partial}^2 + (v_1+v_2) & -(p+q)\bar{\partial}
+ q_{\bar{z}} -2p_{\bar{z}} \\
(p+q)\partial -p_z +2q_z & -(\partial^2+\bar{\partial}^2) -(v_1+v_2)
\end{array}
\right),
$$
where
$$
v_{1\bar{z}} = -2(pq)_z, \ \ \ \
v_{2z} = -2(pq)_{\bar{z}},
$$
the equations (\ref{triple}) takes the form
\begin{equation}
\label{ds2}
\begin{split}
p_t = p_{zz} + p_{\bar{z}\bar{z}} + (v_1+v_2)p,
\\
q_t = -q_{zz}-q_{\bar{z}\bar{z}} -(v_1+v_2)q.
\end{split}
\end{equation}
\end{proposition}

\begin{proposition}
\label{third}
For
$$
A_3 =
\left(
\begin{array}{cc}
\partial^3 + \frac{3}{2}v_1 \partial -3w_1 &
q\bar{\partial}^2 - q_{\bar{z}}\bar{\partial} +q_{\bar{z}\bar{z}} +
\frac{3}{2}v_2 q \\
p\partial^2 - p_z \partial +p_{zz} + \frac{3}{2}v_1 p &
\bar{\partial}^3 + \frac{3}{2} v_2 \bar{\partial} -3w_2
\end{array}
\right),
$$
$$
B_3 =
\left(
\begin{array}{cc}
b_{11} & b_{12} \\
b_{21} & b_{22}
\end{array}\right),
$$
where
$$
b_{11} = -b_{22} =
\bar{\partial}^3  -\partial^3 -
\frac{3}{2}(v_1 \partial - v_2 \bar{\partial})
+3(w_1-w_2),
$$
$$
b_{12} =
-(p+q)\bar{\partial}^2 -\frac{3}{2}(p+q)v_2 -
(3p_{\bar{z}}-q_{\bar{z}})\bar{\partial} -
(3p_{\bar{z}\bar{z}} + q_{\bar{z}\bar{z}}),
$$
$$
b_{21} =
-(p+q)\partial^2 -\frac{3}{2}(p+q)v_1 -
(3q_z-p_z)\partial -
(3q_{zz} + p_{zz})
$$
and
$$
v_{1\bar{z}} = -2(pq)_z, \ \ \ \
v_{2z} = -2(pq)_{\bar{z}},
$$
$$
w_{1\bar{z}} = (pq_z)_z, \ \ \ \
w_{2z} =(qp_{\bar{z}})_{\bar{z}},
$$
the equations (\ref{triple}) takes the form
\footnote{
The equation (\ref{triple}) after a formal substitution of $A$ and $B$
reduces to the system
$$
p_t = p_{zzz} + p_{\bar{z}\bar{z}\bar{z}} +
\frac{3}{2}(v_1 p_z + v_2 p_{\bar{z}}) +
3(w_1-w_2 +\frac{1}{2}v_{1z})p,
$$
$$
q_t = q_{zzz} + q_{\bar{z}\bar{z}\bar{z}} +
\frac{3}{2}(v_1 q_z + v_2 q_{\bar{z}}) -
3(w_1-w_2 -\frac{1}{2}v_{2\bar{z}})q,
$$
}
\begin{equation}
\label{ds3}
\begin{split}
p_t = p_{zzz} + p_{\bar{z}\bar{z}\bar{z}} +
\frac{3}{2}(v_1 p_z + v_2 p_{\bar{z}}) -
3( \partial^{-1}[(qp_{\bar{z}})_{\bar{z}}] +
\bar{\partial}^{-1} [(qp_z)_z]) p,
\\
q_t = q_{zzz} + q_{\bar{z}\bar{z}\bar{z}} +
\frac{3}{2}(v_1 q_z + v_2 q_{\bar{z}}) -
3( \partial^{-1}[(pq_{\bar{z}})_{\bar{z}}] +
\bar{\partial}^{-1} [(pq_z)_z]) q.
\end{split}
\end{equation}
\end{proposition}

The equations (\ref{ds2}) and (\ref{ds3}) are called the Davey--Stewartson
equations. In fact these are the equations DSII$_2$ and DSII$_3$ from the
DSII hierarchy. The equation DSII$_n$ takes the form
(\ref{triple}) where $A_n$ equals
$$
A_n =
\left(
\begin{array}{cc}
(-1)^{n+1} \partial^n  & 0 \\
0 & \bar{\partial}^n
\end{array}
\right)
+ \dots
$$
(here by $\dots$ we denote terms of lower order).

For $n=1$ we have
$$
A_1 =
\left(
\begin{array}{cc}
\partial  & q \\
p & \bar{\partial}
\end{array}
\right), \ \ \
B_1 =
\left(
\begin{array}{cc}
\bar{\partial}- \partial  & -(p+q) \\
-(p+q) & \partial - \bar{\partial}
\end{array}
\right),
$$
and the DSII$_1$ equations are
$$
p_t = p_z+p_{\bar{z}}, \ \ \ \
q_t = q_z + q_{\bar{z}}.
$$

Remark that the DSI hierarchy is a hierarchy of nonlinear equations
obtained from the DSII hierarchy by replacing the variables $z,\bar{z}$
by real-valued variables $x,y$.

Let us consider the reduction of the DSII hierarchy for the case
\begin{equation}
\label{constraint}
p = -u, \ \ \ q = \bar{u}.
\end{equation}

The equation (\ref{ds2}) is not compatible however the substitution
$$
A_2 \to iA_2, \ \ \ B_2 \to iB_2
$$
into (\ref{triple}) gives a reduction of (\ref{ds2}) compatible with
(\ref{constraint}):
\begin{equation}
\label{ds22}
\begin{split}
u_t = i(u_{zz}+u_{\bar{z}\bar{z}} + 2(v+\bar{v})u),
\\
v_{\bar{z}} = (|u|^2)_z.
\end{split}
\end{equation}

The substitution of (\ref{constraint}) into (\ref{ds3}) gives
\begin{equation}
\label{ds32}
\begin{split}
u_t = u_{zzz} + u_{\bar{z}\bar{z}\bar{z}} + 3(v u_z + \bar{v}  u_{\bar{z}}) +
3 (w + w^\prime) u,
\\
v_{\bar{z}} = (|u|^2)_z, \ \ \
w_{\bar{z}} = (\bar{u}u_z)_z, \ \ \
w^\prime_z = (\bar{u}u_{\bar{z}})_{\bar{z}}.
\end{split}
\end{equation}

For brevity we shall call the equations (\ref{ds22}) and (\ref{ds32}) by
the DS$_2$ and DS$_3$ equations respectively.

In difference with (\ref{ds22}) the DS$_3$ equation is
compatible with the constraint $u=\bar{u}$ and for real-valued potentials
it reduces to the modified Novikov--Veselov equation:
\begin{equation}
\label{mnv}
\begin{split}
u_t = u_{zzz} + u_{\bar{z}\bar{z}\bar{z}} + 3(v u_z + \bar{v} u_{\bar{z}}) +
\frac{3}{2} (v_z + \bar{v}_{\bar{z}}) u,
\\
v_{\bar{z}} = (u^2)_z.
\end{split}
\end{equation}

Notice that $A_n$ depends on two functional parameters which are $p$ and $q$
and put
$$
A^+_n = A \ \ \mbox{for $p=-u,q=\bar{u}$}, \ \ \
A^-_n = A \ \ \mbox{for $p=-\bar{u}, q = u$}.
$$

Now let us recall the definition of the DS deformations of a surface
introduced in \cite{K2}.

\begin{proposition}[\cite{K2}]
Let a surface $\Sigma$ be  defined by the formulas (\ref{wr1}) and
(\ref{wr2}) for some $\psi^0,\varphi^0$
and let $U(z,\bar{z},t)$ be a deformation of the potential
described by the equation
(\ref{ds22}) or (\ref{ds32}). Then the formulas
(\ref{wr1}) and
(\ref{wr2}) and the equations
\begin{equation}
\label{dsn}
\begin{split}
\psi_t = iA^+_2 \psi, \ \ \ \varphi_t = - iA^-_2 \varphi,
\\
\psi_t = A^+_3 \psi, \ \ \ \varphi_t = A^-_3 \varphi
\end{split}
\end{equation}
with $\psi_{t=0} = \psi^0, \varphi_{t=0} = \varphi^0$,
define deformations of the surface governed by the equations (\ref{ds22}) and
(\ref{ds32}) respectively.
\end{proposition}

The proof of this proposition is as follows.
Since the deformation of $U=u$ is
described by the equation (\ref{triple}),
the vector functions $\psi$ and $\varphi$ meeting
$\D \psi = 0$ and $\D^\vee\varphi=0$ are deformed via equations (\ref{dsn})
and for any $t$ they meet again the Dirac equations. Therefore by Proposition
\ref{first} they define a surface $\Sigma_t$ via the Weierstrass formulas
(\ref{wr1}) and (\ref{wr2}). Thus we have a deformation $\Sigma_t$ such that
$\Sigma_0 = \Sigma$.

Any equation of the DSII hierarchy defines such a deformation (for $n$ even
we have to substitute $A_n \to i A_n$ to preserve the reduction $p=-\bar{q}$).
We write down only two equations, DS$_2$ and DS$_3$,
because they resemble the main
properties of others and do not involve very large expressions.

For $u=\bar{u}$ such a deformation reduces to the mNV deformation
defined in \cite{K1} and studied in \cite{T1,T12,GL,Langer,BPP}.

\section{The Weierstrass representation}
\label{sec3}

An oriented two-plane in $\R^4$ is defined by a positively-oriented
orthonormal basis
$$
e_1=(e_{1,1},\dots,e_{1,4}), \ \ \
e_2=(e_{2,1},\dots,e_{2,4})
$$
which is defined up to rotations. There is a one-to-one
correspondence
$$
\{(e_1,e_2)\} \leftrightarrow (y_1:y_2:y_3:y_4), \ \ \
y_k = e_{1,k}+i e_{2,k}, \ k=1,2,3,4,
$$
between the moduli space of such planes (which is the Grassmannian
$\widetilde{G}_{4,2}$) and points of the quadric $Q \subset \C P^3$
defined by the equations
$$
y_1^2 + y_2^2 + y_3^2 + y_4^2 = 0.
$$
In terms of another homogeneous coordinates  $y_1^\prime,\dots,y_4^\prime$
such that
$$
y_1 = \frac{i}{2}(y^\prime_1 + y^\prime_2),
\ \ \
y_2 = \frac{1}{2}(y^\prime_1-y^\prime_2),
$$
$$
y_3 = \frac{1}{2}(y^\prime_3 + y^\prime_4),
\ \ \
y_4 = \frac{i}{2}(y^\prime_3-y^\prime_4)
$$
this quadric is written as
$$
y^\prime_1 y^\prime_2 = y^\prime_3 y^\prime_4.
$$
Therefore the correspondence
$$
y^\prime_1 = a_2 b_2,
\ \
y^\prime_2 = a_1 b_1,
\ \
y^\prime_3 = a_2 b_1,
\ \
y^\prime_4 = a_1 b_2
$$
establishes a biholomorphic equivalence
$$
\widetilde{G}_{4,2} = \C P^1 \times \C P^1
$$
where $(a_1:a_2)$ and $(b_1:b_2)$ are homogeneous coordinates on the
copies of $\C P^1$. This mapping $\C P^1 \times \C P^1 \to Q \subset \C P^3$
is called the Segr\'e mapping.

Let $r: W \to \R^4$
be an immersion of a surface with a conformal parameter $z \in W \subset \C$.
The conformality condition reads
$$
\langle r_z, r_z \rangle = \sum_{k=1}^4 (x^k_z)^2 = 0
$$
where $x^k_z = \frac{\partial x^k}{\partial z}$, $k=1,2,3,4$.
The Gauss map takes the form
$$
G: W \to \widetilde{G}_{4,2}, \ \ \
Q \in W \to (x^1_z(Q),:\dots:x^4_z(Q)).
$$
By using the equivalence $\widetilde{G}_{4,2} = \C P^1 \times \C P^1$,
decompose $G$ into two maps
$$
G = (G_\psi,G_\varphi)
$$
where
$$
G_\psi = (\psi_1:\bar{\psi}_2) \in \C P^1, \ \ \
G_\varphi = (\varphi_1 : \bar{\varphi}_2) \in \C P^1
$$
and rewrite the formulas for $x^k_z$ in terms of these maps as follows
\begin{equation}
\label{weier4}
\begin{split}
x^1_z = \frac{i}{2} (\bar{\varphi}_2\bar{\psi}_2 + \varphi_1 \psi_1),
\ \ \ \
x^2_z = \frac{1}{2} (\bar{\varphi}_2\bar{\psi}_2 - \varphi_1 \psi_1),
\\
x^3_z = \frac{1}{2} (\bar{\varphi}_2 \psi_1 + \varphi_1 \bar{\psi}_2),
\ \ \ \
x^4_z = \frac{i}{2} (\bar{\varphi}_2 \psi_1 - \varphi_1 \bar{\psi}_2).
\end{split}
\end{equation}
We have
$$
dx^k = \eta_k, \ \ \ k=1,2,3,4,
$$
where the forms $\eta_k$ take the same shapes as in Proposition \ref{first}.

This decomposition is not unique and functions $\psi$ and $\varphi$
are defined up to gauge transformations
\begin{equation}
\label{gauge1}
\left(
\begin{array}{c}
\psi_1 \\ \psi_2
\end{array}
\right) \to
\left(
\begin{array}{c}
e^f\psi_1 \\ e^{\bar{f}}\psi_2
\end{array}
\right),
\ \ \
\left(
\begin{array}{c}
\varphi_1 \\ \varphi_2
\end{array}
\right) \to
\left(
\begin{array}{c}
e^{-f} \varphi_1 \\ e^{-\bar{f}}\varphi_2
\end{array}
\right)
\end{equation}
By choosing a representative $\psi$ for $G_\psi$ we fix a function
$\varphi$.

The formula (\ref{sec1}) gives exactly a gauge transformation between
different lifts of the Gauss mapping $G=(G_\psi,G_\varphi)$ to non-vanishing
spinors, i.e. to $(\C^2 \setminus \{0\}) \times (\C^2 \setminus \{0\})$,
which we mention in the introduction.

The closedness conditions for the forms $\eta_k$ are
\begin{equation}
\label{closedness}
\left(\bar{\varphi}_2 \psi_1 \right)_{\bar{z}} = \left(\bar{\varphi}_1
\psi_2 \right)_z, \ \ \
\left(\bar{\varphi}_2 \bar{\psi}_2\right)_{\bar{z}} = -
\left(\bar{\varphi}_1 \bar{\psi}_1\right)_z
\end{equation}
and they do not have the form of Dirac equations $\D\psi=0$ and
$\D^\vee\varphi=0$ for arbitrary representatives $(\psi_1,\psi_2)$ and
$(\varphi_1,\varphi_2)$ of the mappings $G_\psi$ and $G_\varphi$.
Such representatives have to be found by solving some differential
equations.

Let us start with the following lift for
$G_\psi = (\psi_1:\bar{\psi}_2)$
which is
correctly defined up to a $\pm 1$ multiple:
$$
s_1 = e^{i\theta} \cos \eta, \ \ \
s_2 = \sin \eta.
$$
We look for a pair of functions $\psi_1,\psi_2$ meeting two conditions:

1) $G_\psi = (\psi_1:\bar{\psi}_2) = (s_1:\bar{s}_2)$;

2) $\D \psi=0$ for some potential $U$.

By the first condition, $\psi$ has the form
$$
\psi_1 = e^g s_1,  \ \ \ \psi_2 = e^{\bar{g}} s_2.
$$
The second condition is written as
$$
\partial(e^{\bar{g}}\sin \eta) + U e^{g+i\theta}\cos \eta = 0,
\ \ \
\bar{\partial}(e^{g+i\theta} \cos \eta) = \bar{U} e^{\bar{g}} \sin \eta.
$$
These equations are rewritten as follows:
$$
U = - \frac{e^{\bar{g}}}{e^{g + i\theta}\cos \eta}
\left( \bar{g}_z \sin \eta +
\eta_z \cos \eta \right),
$$
$$
U = \frac{e^{\bar{g}}}{e^{g + i\theta}\sin \eta}
\left( \bar{g}_z \cos \eta - i \theta_z
\cos \eta - \eta_z \sin \eta \right),
$$
which imply
$$
g_{\bar{z}} + i \theta_{\bar{z}} \cos^2 \eta = 0,
$$
$$
U = - e^{\bar{g}-g-i\theta}\left( i \theta_z \sin \eta \cos \eta +
\eta_z \right).
$$

The $\bar{\partial}$-problem for $g$ is solved by the well-known means and
its solution is defined up to holomorphic functions.
Therefore the potential $U$ is defined up to a multiplication by
$e^{\bar{h}-h}$ where $h$ is an arbitrary holomorphic function.

It is derived by straightforward computations that $\D\psi=0$ implies that
the condition (\ref{closedness}) takes the form
of the Dirac equation $\D^\vee\varphi = 0$.

Thus the following theorem is derived.

\begin{theorem}
\label{theorem4}
Let $r:W \to \R^4$ be an immersed surface with a conformal parameter
$z$ and let
$G_\psi = (e^{i\theta}\cos \eta: \sin \eta)$ be one of the components
of its Gauss map.

There exists another representative $\psi$ of this mapping $G_\psi =
(\psi_1:\bar{\psi}_2)$ such that it meets the Dirac equation
$$
\D \psi = 0
$$
with some potential $U$.

A vector function $\psi = (e^{g+i\theta}\cos \eta,e^{\bar{g}}\sin \eta)$
is defined from the equation
\begin{equation}
\label{maineq}
g_{\bar{z}} = -i\theta_{\bar{z}} \cos^2 \eta,
\end{equation}
up to holomorphic functions and the corresponding potential $U$
is defined up  by the formula
$$
U =  -e^{\bar{g}-g-i\theta}(i\theta_z \sin \eta \cos \eta +\eta_z)
$$
up to multiplications by $e^{\bar{h}-h}$ where $h$ is an arbitrary
holomorphic function.

Given the function $\psi$, a function $\varphi$
which represents another component $G_\varphi$ of the Gauss map
meets the equation
$$
\D^\vee \varphi = 0.
$$

Different representations (lifts) of the Gauss mapping $G$ of the surface
$W$ are related by gauge transformations of the form
\begin{equation}
\label{gauge2}
\begin{split}
\left(
\begin{array}{c}
\psi_1 \\ \psi_2
\end{array}
\right) \to
\psi^\prime = \left(
\begin{array}{c}
e^h\psi_1 \\ e^{\bar{h}}\psi_2
\end{array}
\right),
\ \ \
\left(
\begin{array}{c}
\varphi_1 \\ \varphi_2
\end{array}
\right) \to
\varphi^\prime =
\left(
\begin{array}{c}
e^{-h} \varphi_1 \\ e^{-\bar{h}}\varphi_2
\end{array}
\right),
\\
U \to U^\prime = e^{\bar{h}-h}U, \
\end{split}
\end{equation}
where $h$ is an arbitrary holomorphic function on $W$.
 \end{theorem}

Since we decompose the Gauss map into
two components $G_\psi$ and $G_\varphi$
with functions $\psi$ and $\varphi$ meeting the Dirac equations, we obtain
the surface by integrating the differentials $dx^k$ given by (\ref{weier4}).
The formulas (\ref{metric}) and (\ref{potential}) for the metric and
the potential are obtained by simple computations. We have

\begin{corollary}
Every oriented surface in $\R^4$ admits a Weierstrass
representation given by Proposition \ref{first}.
\end{corollary}

However the non-uniqueness of such a representation leads to
the following conclusion.

\begin{corollary}
For any surface $r: W \to \R^4$ where $W$ is an open subset of $\C$
and any subdomain $V \subset W$
such that $V \neq W$ there is a Weierstrass representation of
$r\vert_V:V \to \R^4$
which is not analytically continued (i.e. expanded) onto $W$.
\end{corollary}

{\sc Proof.} For that take a Weierstrass representation of $W$, restrict it
onto $V$ and take a gauge equivalent representation of $V$
corresponding to a transformation (\ref{gauge2}) where $h: V \to \C$ is a
holomorphic function which is not analytically continued onto $W$.
Then this representation of $V$ is not expanded onto $W$.
This proves the corollary.

The following theorem is clear.

\begin{theorem}
Given a Weierstrass representation of an immersed
closed oriented surface $\Sigma$ into $\R^4$, the corresponding functions
$\psi$ and $\varphi$ are sections of the $\C^2$-bundles
$E$ and $E^\vee$ over $\Sigma$ which are as follows:

1) $E$ and $E^\vee$ split into sums of pair-wise conjugate line bundles
$$
E = E_0 \oplus \bar{E}_0, \ \ \ \ E^\vee = E^\vee_0 \oplus \bar{E}^\vee_0
$$
such that $\psi_1$ and $\bar{\psi}_2$ are sections of $E_0$ and
$\varphi_1$ and $\bar{\varphi}_2$ are sections of $E^\vee_0$;

2) the pairing of sections of $E_0$ and $E^\vee_0$ is a $(1,0)$ form
on $\Sigma$: if
$$
\alpha \in \Gamma(E_0), \ \ \ \beta \in \Gamma(E^\vee_0),
$$
then
$$
\alpha \beta dz
$$
is a correctly defined $1$-form on $\Sigma$;

3) the Dirac equation $\D \psi =0$ implies that $U$ is a section of
the same line bundle $E_U$ as
$$
\frac{\partial \gamma}{\alpha} \in \Gamma(E_U) \ \ \
\mbox{for $\alpha \in \Gamma(E_0),\ \  \gamma \in \Gamma(\bar{E}_0)$}
$$
and the quantity $U\bar{U} dz \wedge d\bar{z}$ is a correctly defined
$(1,1)$-form on $\Sigma$ whose integral equals
$$
\int_{\Sigma} U\bar{U} dz \wedge d\bar{z} = - \frac{i}{2} {\cal W}(\Sigma)
$$
where ${\cal W}(\Sigma) = \int_{\Sigma} |{\bf H}|^2 d\mu$ is the Willmore
functional of $\Sigma$.
\end{theorem}

Notice that for surfaces in $\R^3$ functions $\psi$ are sections of
spinor bundles \cite{T1} and the gauge transformation (\ref{gauge2}) shows
that for surfaces $\R^4$ this is not necessarily a case.

It is derived from Proposition \ref{first} that,
given a Riemann surface $\Sigma$, such bundles
$E,E^\vee$, and $E_U$ and solutions $\psi$ and $\varphi$ to the equations
$\D \psi =0$ and $\D^\vee \varphi =0$, one may construct an immersion of
the universal covering of $\Sigma$ into $\R^4$.
The Gauss mapping of this immersion descends through $\Sigma$.
Let us give a criterion for converting such an immersion into an immersion of
the surface $\Sigma_0$.

\begin{proposition}
Let $\Sigma$ be an oriented closed surface, let
$\widetilde{\Sigma}$ be its universal covering  and let $(\psi,\varphi)$
define an immersion of $\widetilde{\Sigma}$
into $\R^4$ via (\ref{wr1}) and (\ref{wr2}).
Then such an immersion converts
into an immersion of $\Sigma$ if and only if
\begin{equation}
\label{period4}
\int_{\Sigma} \bar{\psi}_1\bar{\varphi}_1 d\bar{z} \wedge \omega =
\int_{\Sigma} \bar{\psi}_1 \varphi_2 d\bar{z} \wedge \omega =
\int_{\Sigma} \psi_2 \bar{\varphi}_1 d\bar{z} \wedge \omega =
\int_{\Sigma} \psi_2 \varphi_2 d\bar{z} \wedge \omega =0
\end{equation}
for any holomorphic differential $\omega$ on $\Sigma$.
\end{proposition}

{\sc Proof.} Let $g \geq 1$ be the genus of $\Sigma$. Then there is a basis
$\alpha_1,\dots,\alpha_g$, $\beta_1$,$\dots$,$\beta_g$
of $1$-cycles on $\Sigma$
such that $\Sigma = \widetilde{\Sigma}/\Gamma$ and the fundamental domain
for the action of $\Gamma = \pi_1(\Sigma)$ is
a domain $\Omega$ on $\widetilde{\Sigma}$ whose boundary has
the form
$$
\partial \Omega = \alpha_1 \beta_1 \alpha_1^{-1} \beta_1^{-1}
\dots \alpha_g \beta_g \alpha^{-1} \beta_g^{-1}.
$$
Denote by $dx^1,\dots,dx^4$ the closed $1$-forms on $\widetilde{\Sigma}$
induced by the immersion into $\R^4$ and denote by $V^1,\dots,V^4$
the period vectors
$$
V^k = \left(\int_{\alpha_1}dx^k,\dots,\int_{\alpha_g}dx^k,
\int_{\beta_k}dx^k,\dots,\int_{\beta_g}dx^k\right), \ \ \
k=1,2,3,4.
$$
An immersion of $\widetilde{\Sigma}$ converts into an immersion of
$\Sigma$ if and only if
$$
V^1 = V^2 = V^3 = V^4 = 0.
$$

Given a holomorphic form $\omega$ on $\Sigma$, pull it back onto the
universal covering $\widetilde{\Sigma}$ and compute the integral
$$
\int_{\partial \Omega} x^k \omega =
\sum_{j=1}^g \left(\int_{\alpha_j} \omega \int_{\beta_j} dx^k -
\int_{\alpha_j}dx^k \int_{\beta_j} \omega \right) =
$$
$$
=
\sum_{j=1}^g \left(V^k_{j+g} \int_{\alpha_j} \omega - V^k_j
\int_{\beta_j} \omega\right)
$$
which is by the Stokes theorem equals
$$
\int_{\partial \Omega} x^k \omega =
\int_\Omega x^k_{\bar{z}} d\bar{z}\wedge \omega.
$$

A Riemann surface $\Sigma$ has
a basis $\omega_1,\dots,\omega_g$
for holomorphic differentials normalized by the condition
$$
\int_{\alpha_j} \omega_k = \delta_{ij}.
$$
In this event the $\beta$-periods matrix
$$
B_{jk} = \int_{\beta_j} \omega_k
$$
is symmetric with positive imaginary part: $\Im B >0$.
This implies that the conditions
$$
\sum_{j=1}^g \left(V^k_{j+g} \int_{\alpha_j} \omega_l - V^k_j
\int_{\beta_j} \omega_l\right) = 0, \ \ \ l=1,\dots,g,
$$
for a vector $V^k$ with real entries are satisfied if and only if
$V^k=0$. Therefore an immersion of $\widetilde{\Sigma}$ converts
into an immersion of a closed surface $\Sigma$ if and only if
$$
\int_\Omega x^k_{\bar{z}} d\bar{z} \wedge \omega = 0
$$
for any $k=1,\dots,4$ and any holomorphic differential $\omega$ on $\Sigma$.
It follows from (\ref{wr1}) that that is equivalent to
the equalities (\ref{period4}). This proves the proposition.

\section{Deformations of tori via the Davey--Ste\-wart\-son equations}
\label{sec4}

Let us look what Theorem \ref{theorem4} gives us for tori.

\begin{theorem}
\label{torus}
Let $\Sigma$ be a torus in $\R^4$ which is conformally equivalent to
$\C/\Lambda$ and $z$ is a conformal parameter.

Then there are vector functions $\psi,\varphi$ and a function $U$
such that

1) $\psi$ and $\varphi$ give a Weierstrass representation of
$\Sigma$;

2) the potential $U$ of this representation is
$\Lambda$-periodic;

3) such functions $\psi$, $\varphi$, and $U$ are defined up to
gauge transformations
\begin{equation}
\label{gauge3}
\left(
\begin{array}{c}
\psi_1 \\ \psi_2
\end{array}
\right)
\to
\left(
\begin{array}{c}
e^{h} \psi_1 \\ e^{\bar{h}} \psi_2
\end{array}
\right),
\ \
\left(
\begin{array}{c}
\varphi_1 \\ \varphi_2
\end{array}
\right)
\to
\left(
\begin{array}{c}
e^{-h} \varphi_1 \\ e^{-\bar{h}}\varphi_2
\end{array}
\right),
\ \
U \to U e^{\bar{h} -h}U
\end{equation}
where
$$
h(z) = a + bz, \ \ \Im (b\gamma) \in \pi \Z \ \
\mbox{for all $\gamma \in \Lambda$}.
$$
Therefore such representations are parameterized by a $\Z^2$ lattice
formed by admissible values of $b$.
\end{theorem}

{\sc Proof.} By Theorem \ref{theorem4}, given a Weierstrass
representation of a torus (with a fixed conformal parameter)
the Dirac equations are satisfied if and only if
$$
U = - e^{\bar{g}-g-i\theta}\left( i \theta_z \sin \eta \cos \eta +
\eta_z \right), \ \
\psi_1 = e^{i\theta+g}\cos \eta, \ \ \psi_2 = e^{\bar{g}}\sin\eta
$$
where
$g$ meets the equation (\ref{maineq}).
The component of the Gauss mapping
$G_\psi =(e^{i\theta}\cos\eta:\sin \eta)$ is $\Lambda$-periodic and
$U$ takes the form
$U = e^{\bar{g}-g} U_0$
where
$$
U_0 = i \theta_z \sin \eta \cos \eta +
\eta_z
$$
is periodic with respect to $\Lambda$.
By (\ref{maineq}), $g_{\bar{z}}$ is periodic
and  a general solution to (\ref{maineq}) takes the form
$$
g = h(z) + c\bar{z} + f(z,\bar{z})
$$
where $h(z)$ is an arbitrary holomorphic function, $f$ is some periodic
periodic function, and
$$
c = - \int_M i \theta_{\bar{z}} \cos^2 \eta dx dy.
$$
Hence $U$ is a $\Lambda$-periodic function if and only if
$g = a + bz + c \bar{z}$ (modulo periodic functions) and
$(\bar{g}(z) - g(z)) \in 2 \pi \Z$ for $z \in \Lambda$.
The latter conditions are easily resolved and $b$ is defined up to
$b^\prime$ such that $(\overline{b^\prime \gamma} -b^\prime \gamma)
\in 2\pi \Z$ for all $\gamma \in \Lambda$.
This proves the theorem.

Now let us look for the DS deformations of tori.
Exactly in this case in difference with high genera compact surfaces
the constraints for defining additional potentials $v,a,b$ could be globally
resolved. As in \cite{T1} we are interested in two problems:

\begin{itemize}
\item
when a deformation from the DS hierarchy deforms a torus into tori?

\item
when the ``Willmore functional'' $\int_{\Sigma} u\bar{u} dz \wedge d\bar{z}$
is preserved by such a deformation.
\end{itemize}

For that we have to resolve the constraint equations for $v,w$, and $w^\prime$
carefully choosing specific solutions.

We consider only the  DS$_2$ and DS$_3$ equations
since we do not know until recently a general recursion procedure
for writing down explicit formulas for higher equations.
However it is a general point in the soliton theory that the first nontrivial
equations in the hierarchy usually resemble the main properties of
higher equations.

Since we derive for tori that there is a representation with a double-periodic
potential we can look for solutions of the DS equations with such an initial
data. In addition we have to define such resolutions of constraints (i.e.
the additional potentials $v,a,b$) to obtain the equations with
double-periodic coefficients and, hence, double-periodic solutions.

We do not discuss the existence of a solution assuming that it exists
which, in particular, for short times follows from the Cauchy--Kovalevskaya
theorem.

{\sc The DS$_2$ deformation.}

We have
$$
\psi_t = iA_2^+ \psi, \ \ \varphi_t = -iA_2^- \varphi, \ \ v_{\bar{z}} =
(|u|^2)_z.
$$

When we work with compact surfaces we have to resolve the constraint equation
for $v$ globally. Moreover for tori we would like to
save the periodicity of the integrands
of the Weierstrass representation, i.e. terms which are
of  the form $\psi_2 \varphi_2$ or similar to it, we have to have a periodic
potential $v$.

This is quite easy: the constraint for $v$ is uniquely resolved
by an inversion of the $\bar{\partial}$-operator on a torus assuming that
$\int v dz \wedge d\bar{z} =0$. The reason is the same as for the mNV
deformation (see \cite{T1}): since the right hand-side $(|u|^2)_z$ of
the constraint equation is a derivative of a periodic function its Fourier
decomposition does not contain a non-zero
term and the $\bar{\partial}$-operator
is inverted term by term of the Fourier decomposition.
Thus we have:
\begin{equation}
\label{cons1}
v = \bar{\partial}^{-1}\partial (|u|^2), \ \ \
\int_{\Sigma} v dz \wedge d\bar{z} = 0.
\end{equation}
Here and in the sequel we identify $\Sigma$ with a fundamental domain for the
lattice $\Lambda$ such that the torus is conformally equivalent to
$\C/\Lambda$.

\begin{theorem}
\label{2tori}
The DS$_2$ equation
$$
u_t = i(u_{zz}+u_{\bar{z}\bar{z}} + 2(v+\bar{v})u)
$$
where $v$ is defined by (\ref{cons1})
induces a deformation of tori (with a fixed periodic potential of their
Weierstrass representations) into tori.
\end{theorem}

{\sc Proof.}
We have to prove that the closedness conditions (\ref{period4})
are preserved. We show that only for one of them because for others the proofs
are basically the same. There is one dimensional family of holomorphic
differential on a torus generated by $dz$. We have to prove that
$$
J = \int_{\Sigma} \psi_2 \varphi_2 dz \wedge d\bar{z} = 0
$$
implies that $J_t=0$.

We have
$$
\psi_{2t} = i[(u\partial -u_z) \psi_1 + (\bar{\partial}^2 + \bar{v})\psi_2],
\ \
\varphi_{2t} = -i[(\bar{u} \partial - \bar{u}_z)\varphi_1 +
(\bar{\partial}^2 + \bar{v})\varphi_2].
$$
Substituting that into
$$
J_t = \int_{\Sigma} (\psi_{2t}\varphi_2 +
\psi_2\varphi_{2t}) dz \wedge d\bar{z}
$$
(notice that the integrand is correctly defined as a function on $\Sigma$, i.e.
it is double-periodic, although $\psi_2$ and $\varphi_2$ are not periodic)
we obtain
$$
J_t =
i\int_{\Sigma}
[(u \psi_{1z} -u_z \psi_1)\varphi_2 + (\psi_{2\bar{z}\bar{z}} + \bar{v}\psi_2)
\varphi_2) -
$$
$$
(\bar{u}\varphi_{1z} - \bar{u}_z \varphi_1)\psi_2 -
(\varphi_{2\bar{z}{\bar z}} + \bar{v}\varphi_2) \psi_2] dz \wedge d\bar{z}.
$$
Integrating by parts some of terms
and canceling terms
repeated with different signs we derive
$$
J_t =
i\int_{\Sigma}
[(u \psi_{1z} -u_z \psi_1)\varphi_2 -
(\bar{u}\varphi_{1z} - \bar{u}_z \varphi_1)\psi_2]
dz \wedge d\bar{z}.
$$
Now an integration by parts implies that
$$
J_t = i \int_{\Sigma}
[u \psi_{1z}\varphi_2 + u(\psi_1 \varphi_2)_z - \bar{u}\varphi_{1z}\psi_2
- \bar{u}(\psi_2 \varphi_1)_z]dz \wedge d\bar{z} =
$$
$$
= \int_{\Sigma}
[2(u \varphi_2) \psi_{1z} + (u\psi_1) \varphi_{2z} -
2(\bar{u}\psi_2)\varphi_{1z} - (\bar{u}\varphi_1)\psi_{2z}]dz \wedge d\bar{z}.
$$
By using the Dirac equations replace the terms in brackets to
exclude the potentials  $u$ and $\bar{u}$ from the integrand:
$$
J_t = 2i\int_{\Sigma} (\varphi_{1\bar{z}}\psi_{1z} - \varphi_{1z}
\psi_{1\bar{z}}) dz \wedge d\bar{z}.
$$
An integration by parts implies $J_t = 0$ which proves the theorem.

Remark that we never use in the proof that $J=0$ for the initial torus.
Therefore we proved more:

{\sl the DS$_2$ deformation preserves the translational periods of
surfaces with double-periodic Gauss mapping.}

Since we proved that tori are preserved it is reasonable to speak about the
conservation laws for the DS deformations. We have

\begin{theorem}
\label{2willmore}
The DS$_2$ deformation of tori preserves the Willmore functional.
\end{theorem}

{\sc Proof.} We have
$$
\frac{d}{dt} \int_{\Sigma} |u|^2 dz \wedge d\bar{z} =
\int_{\Sigma}(u_t \bar{u} + u \bar{u}_t) dz \wedge d\bar{z} =
$$
$$
=
 i \int_{\Sigma} [\bar{u}(u_{zz} + u_{\bar{z}\bar{z}} + 2(v+ \bar{v})u) -
u(\bar{u}_{zz} + \bar{u}_{\bar{z}\bar{z}} + 2(v+\bar{v})\bar{u})]
dz \wedge d\bar{z} =
$$
$$
= i \int_{\Sigma} (\bar{u} u_{zz} - u \bar{u}_{zz} + \bar{u} u_{\bar{z}\bar{z}}
- u \bar{u}_{\bar{z}\bar{z}}) dz \wedge d\bar{z}.
$$
By integrating by parts the last integral, we easily derive that
$$
\frac{d}{dt} \int_{\Sigma} |u|^2 dz \wedge d\bar{z} = 0
$$
which proves the theorem.

{\sc Remark.} By Proposition \ref{second}, there are two potentials $v_1$
and $v_2$ which are coming into the equation. We choose an additional
constraint $v = v_1 = \bar{v}_2$. Without this constraint Theorem \ref{2tori}
does not hold as one can see from its proof.

In fact we assumed more, that
$\int v_{\Sigma} dz \wedge d\bar{z}$ vanishes. This was not important:
if we put $v \to v + f(t)$ then the deformation reduces to a composition of
the initial one and
some evolution tangent to a torus, i.e. this results in adding
a diffeomorphism of a torus which does not affect the geometric picture.

{\sc The DS$_3$ deformation.}

For this deformation we have
$$
\psi_t = A^+_3 \psi, \ \ \ \varphi_t = A^-_3 \varphi
$$
and many constraints which we have to resolve.
As in the case of the DS$_2$ deformation we put
\begin{equation}
\label{ds3c1}
v_{\bar z} = (|u|^2)_z, \ \ \ \int_{\Sigma} v dz \wedge d\bar{z} = 0,
\end{equation}
and in addition we choose $w$ and $w^\prime$ as follows:
\begin{equation}
\label{ds3c2}
w = \partial \bar{\partial}^{-1} (\bar{u}u_z), \ \ \
w^\prime = \bar{\partial} \partial^{-1} (\bar{u}u_{\bar{z}}).
\end{equation}
These functions satisfy the constraint equations (\ref{ds32})
however they are very specific particular solutions to them.

Moreover we have to resolve the constraint equations for $w_1$ and $w_2$
which are different for $A^+$ and $A^-$. We put
\begin{equation}
\label{ds3c3}
\begin{split}
w_1^+ = w -v_z, \ \ \ w_2^+ = -w^\prime,
\\
w_1^- = - w, \ \ \ w_2^- = w^\prime - \bar{v}_{\bar{z}},
\\
v_1^\pm = 2v, \ \ \ v_2^\pm = 2\bar{v}.
\end{split}
\end{equation}
The reasonings for these choices we shall explain later.

\begin{theorem}
\label{3tori}
The DS$_3$ equation
$$
u_t = u_{zzz}+u_{\bar{z}\bar{z}\bar{z}} + 3(v u_z +
\bar{v} u_{\bar{z}}) + 3(w+w^\prime) u
$$
induces a deformation of tori (with fixed periodic potential of their
Weierstrass representations) into tori.
\end{theorem}

{\sc Proof.} We again demonstrate that only for one of the closedness
conditions. Let us take the same as in the proof of Theorem \ref{2tori}.
We have
$$
\psi_{2t} = (-u\partial^2 + u_z \partial -u_{zz} - 3vu)\psi_1 +
(\bar{\partial}^3 + 3 \bar{v} \bar{\partial} + 3 w^\prime)\psi_2,
$$
$$
\varphi_{2t} = (-\bar{u} \partial^2 + \bar{u}_z \partial - \bar{u}_{zz} -
3v\bar{u})\varphi_1 +
(\bar{\partial}^3 + 3\bar{v}\bar{\partial} - 3(w^\prime -\bar{v}_{\bar{z}}))
\varphi_2.
$$
Substitute that into
$$
J_t = \int_{\Sigma} (\psi_{2t} \varphi_2 + \psi_2 \varphi_{2t})
dz \wedge d\bar{z}.
$$
An integration by parts shows that
$$
\int_{\Sigma}
[(\psi_{2\bar{z}\bar{z}\bar{z}} + 3 \bar{v}\psi_{2\bar{z}})\varphi_2 +
(\varphi_{2\bar{z}\bar{z}\bar{z}} + 3 \bar{v}\varphi_{2\bar{z}} +
3\bar{v}_{\bar{z}} \varphi_2)\psi_2] dz \wedge d\bar{z} = 0
$$
and we are left to prove that
$$
\int_{\Sigma}
[
(-u\psi_{1zz} + u_z \psi_{1z} -(u_{zz}+ 3vu)\psi_1)\varphi_2
+
$$
$$
+
(-\bar{u} \varphi_{1zz} + \bar{u}_z \varphi_{1z} -
(\bar{u}_{zz} + 3v\bar{u})\varphi_1)\psi_2
]
dz \wedge d\bar{z} = 0.
$$
Let us rewrite the left-hand side of this formula as
$$
J_t = \int_{\Sigma} [ (-(u\psi_1)_{zz} + 3 u_z \psi_{1z} - 3vu \psi_1)
\varphi_2 +
$$
$$
(-(\bar{u}\varphi_1)_{zz} +
3\bar{u}_z \varphi_{1z} -3v\bar{u}\varphi_1)\psi_2)]
dz \wedge d\bar{z}
$$
which, by the Dirac equations $\D\psi = \D^\vee\varphi=0$, equals to
$$
J_t = \int_{\Sigma} [(\psi_{2zzz} + 3u_z \psi_{1z})\varphi_2  -3v \psi_1
\varphi_{1\bar{z}} +
$$
$$
(\varphi_{2zzz} + 3\bar{u}_z \varphi_{1z})\psi_2  - 3 v \psi_{1\bar{z}}
\varphi_1]
dz \wedge d\bar{z} =
$$
$$
= \int_{\Sigma} (\psi_{2zzz}\varphi_2+\psi_2\varphi_{2zzz}) dz \wedge d\bar{z}
+ 3 \int_{\Sigma} (u_z \psi_{1z}\varphi_2 + \bar{u}_z \psi_2\varphi_{1z}
- v (\psi_1\varphi_1)_{\bar{z}}) dz \wedge d\bar{z}.
$$
An integration by parts shows that the first summand vanishes and, by
(\ref{ds3c1}), we have
$$
J_t =
3 \int_{\Sigma} (u_z \psi_{1z}\varphi_2 + \bar{u}_z \psi_2 \varphi_{1z} +
(|u|^2)_z \psi_1\varphi_1) dz \wedge d\bar{z} =
$$
$$
3 \int_{\Sigma} (u_z \psi_{1z}\varphi_2 + \bar{u}_z \psi_2 \varphi_{1z} +
u_z \bar{u} \psi_1\varphi_1 + \bar{u}_z u \psi_1 \varphi_1) dz \wedge d\bar{z}
=
$$
$$
3 \int_{\Sigma} [u_z(\psi_{1z}\varphi_2 - \psi_1\varphi_{2z}) +
\bar{u}_z (\varphi_{1z}\psi_2 - \varphi_1\psi_{2z})] dz \wedge d\bar{z} =
$$
$$
3\int_{\Sigma}
[u(\psi_1\varphi_{2zz} -\psi_{1zz}\varphi_2) + \bar{u}(\varphi_1\psi_{2zz}
-\varphi_{1zz}\psi_2)] dz \wedge d\bar{z} =
$$
$$
3\int_{\Sigma}( -\psi_{2z}\varphi_{2zz} - \psi_{1zz}\varphi_{1\bar{z}}
- \varphi_{2z} \psi_{2zz} - \psi_{1\bar{z}}\varphi_{1zz}) dz \wedge d\bar{z}.
$$
Now by integrating by parts we easily derive from the last formula that
$$
J_t = 0.
$$
Analogous reasonings show that other closedness conditions (\ref{closedness})
are also preserved by the flow. This proves theorem.

As in the case of the DS$_2$ flow we actually show that

{\sl the DS$_3$ deformation preserves the translational periods of surfaces
with double-periodic Gauss mapping.}

As in the case of the DS$_2$ deformation the Willmore functional is preserved,
i.e. the following theorem holds.

\begin{theorem}
\label{3willmore}
The DS$_3$ deformation with the additional potentials of the form
(\ref{ds3c1}-\ref{ds3c3}) preserves the Willmore functional.
\end{theorem}

{\sc Proof.} We have
$$
\frac{d}{dt} \int_{\Sigma} |u|^2 dz \wedge d\bar{z} =
\int_{\Sigma}(u_t \bar{u} + u \bar{u}_t) dz \wedge d\bar{z} =
$$
$$
\int_{\Sigma}[
(u_{zzz}+u_{\bar{z}\bar{z}\bar{z}} + 3(v u_z +
\bar{v} u_{\bar{z}}) + 3(w+w^\prime) u)\bar{u} +
$$
$$
+
u(\bar{u}_{zzz}+ \bar{u}_{\bar{z}\bar{z}\bar{z}} +
3(\bar{v} \bar{u}_{\bar{z}} +
v \bar{u}_z) + 3(\bar{w}+\bar{w}^\prime) \bar{u})] dz \wedge d\bar{z} =
$$
$$
\int_{\Sigma} [(u_{zzz}\bar{u}+u\bar{u}_{zzz}) +
(u_{\bar{z}\bar{z}\bar{z}}\bar{u} + u \bar{u}_{\bar{z}\bar{z}\bar{z}})]
dz \wedge d\bar{z} +
$$
$$
+
3 \int_{\Sigma} (v(|u|^2)_z + \bar{v}(|u|^2)_{\bar{z}})dz \wedge d\bar{z}
+
3 \int_{\Sigma} (w + \bar{w} + w^\prime + \bar{w}^\prime)|u|^2
dz \wedge d\bar{z}.
$$
An integration by parts shows that the first integral vanishes, and, by
(\ref{ds3c1}), the second integral equals to
$$
3 \int_{\Sigma} (v v_{\bar{z}} + \bar{v}\bar{v}_z) dz\wedge d\bar{z} = 0
$$
(here we use that the function $v$ is double-periodic). We are left to prove
that
$$
\int_{\Sigma} (w + \bar{w} + w^\prime + \bar{w}^\prime)|u|^2
dz \wedge d\bar{z} =0.
$$
We have
$$
w + \bar{w}^\prime = v_z, \ \ \
w^\prime+\bar{w} = \bar{v}_{\bar{z}}.
$$
Therefore the investigated integral is rewritten as
$$
\int_{\Sigma} (v_z + \bar{v}_{\bar{z}})|u|^2 dz \wedge d \bar{z} =
- \int_{\Sigma}(v(|u|^2)_z + \bar{v}(|u|^2)_{\bar{z}}) dz \wedge d\bar{z} =
$$
$$
-\int_{\Sigma} (vv_{\bar{z}} + \bar{v}\bar{v}_z) dz \wedge d\bar{z} = 0.
$$
This proves the theorem.

{\sc Remark.} As for the DS$_2$ deformation we choose the additional
potentials carefully to make the DS$_3$ deformation geometric. For
general potentials meeting the constraint equations it is not the case.
For instance, we may add any constants $c$ and $c^\prime$
to $w$ and $w^\prime$ but the proof of
Theorem \ref{3willmore} shows that if
$k = c+c^\prime+\bar{c}+\bar{c}^\prime \neq 0$ then the Willmore functional
is evolved as
$$
{\cal W}_t = 3k{\cal W}.
$$
and is not preserved.

The most interesting feature of the DS deformations is that they are defined
only for surfaces with fixed potentials of their Weierstrass representations.
Indeed, for a torus we may take another gauge-equivalent potential
\begin{equation}
\label{gp}
u \to u^\prime = e^{az-\bar{a}\bar{z}}u
\end{equation}
and apply the DS deformation for a torus with the potential.
In this case the deformation would be completely different geometrically.
It is noticeable from the deformation of
$|u|^2$ which is, by Proposition \ref{first}, is a geometric quantity.

Let us demonstrate that for the DS$_3$ deformation. The additional
potentials defined by (\ref{ds3c1}-\ref{ds3c3}) are the same as for $u$ but
at $t=0$ the deformation of $|u^\prime|^2$ is
different from the deformation of $|u|^2$ and it is as follows:
$$
\frac{d |u^\prime|^2}{dt} = \frac{d |u|^2}{dt} + 6 \, \mathrm{Re}\,
[a^2(u_z\bar{u} + u\bar{u}_z) + a(u_{zz}\bar{u} - u\bar{u}_{zz})].
$$
Although the first additional term is simple and equals to
$$
3 \left(a^2 \frac{\partial |u|^2}{\partial z}
+ \bar{a}^2 \frac{\partial |u|^2}{\partial \bar{z}}\right),
$$
i.e. could come from one-parametric diffeomorphism group of the surface,
the second term involves the second derivatives and
does not have such a form.

We conclude that

\begin{itemize}
\item
{\sl the DS deformations are correctly defined only for surfaces with fixed
potentials of their Weierstrass representations and for different choices of
the potentials such deformations are geometrically different.}
\end{itemize}

By Theorem \ref{torus}, for tori such deformations are parameterized by a
$\Z^2$ lattice.

If we shall speak on local deformations then the gauge group is much larger
(a gauge transformation is determined by a holomorphic function) and
local deformations would be very different for different choices of
gauge-equivalent potentials.

\section{The spectral curve}

It is reasonable to define the spectral curve for a torus in $\R^4$
as the spectral curve of a double-periodic
operator $\D$ coming in its Weierstrass representation.

Let us recall the definition of the spectral curve of a
double-periodic Dirac operator $\D$ (\cite{T2}).

For that consider all  formal solutions to the equation
$$
\D \psi = 0
$$
meeting in addition the following periodicity conditions:
$$
\psi(z+\gamma_j) =
e^{2\pi i (k_1 \Re \gamma_j + k_2 \Im \gamma_j)}
\psi(z) = \mu(\gamma_j)\psi(z)
\ \ j=1,2.
$$
where $z \in \C$ and $\gamma_1,\gamma_2$ generate the periods lattice
$\Lambda \subset \C$.
Such solutions $\psi$ are called Floquet eigenfunctions
(on the zero level of energy),
the quantities $k_1,k_2$ are called the quasimomenta of $\psi$,
and $(\mu_1,\mu_2) = (\mu(\gamma_1)$, $\mu(\gamma_2)$ are the multipliers of
$\psi$.

The quasimomenta satisfy some analytic relation
(called in solid physics the dispersion relation):
$$
P(k_1,k_2) = 0
$$
which defines a complex curve $Q_0$
in $\C^2$ invariant with respect to translations
$$
k \to k + \gamma^\ast, \ \ \ \gamma^\ast=(\gamma_1^\ast,\gamma_2^\ast)
\in \Lambda^\ast,
$$
where $\Lambda^\ast \subset \C = \R^2 \subset \C^2$ is the dual lattice to
$\Lambda$.

We say that the complex curve $\Gamma = Q_0/\Lambda^\ast$ is
the spectral curve
of $\D$ (on the zero energy level). This definition originates in
the definition of such a curve for a two-dimensional Schr\"odinger operator
\cite{DKN}. The mapping ${\cal M}: \Gamma \to \C^2$ formed
by the multipliers
${\cal M}=(\mu_1,\mu_2)$  is called the multiplier mapping.

The spectral genus of a torus is defined as the geometric genus of the
normalization of $\Gamma$.

For tori in $\R^3$ it appears that such a curve together with ${\cal M}$
contains an important information about the conformal geometry of a torus.
Our conjecture confirmed in \cite{GS} reads that the pair $(\Gamma,{\cal M})$
is preserved by conformal transformations of the ambient space $\R^3$ which
map the torus into $\R^3$. The discussion of other properties of the spectral
curve can be found in \cite{T3}.

For tori in $\R^4$ the situation is slightly different:
the curve $\Gamma$ is defined up to biholomorphic equivalences
however the multiplier mappings depend on the choice of a potential:
the gauge transformation (\ref{gp}) acts on $\psi$ and
${\cal M}$ as follows:
$$
\psi \to e^{-az}\psi,
$$
$$
(\mu_1,\mu_2) \to (e^{-a\gamma_1}\mu_1,e^{-a\gamma_2}\mu_2).
$$

This is rather reasonable. In the fundamental paper \cite{N} by Novikov
the spectral curve of an operator with a potential deformed via some
soliton equation was considered as a conservation law itself for this
equation. Since we show in \S \ref{sec4} that there are infinitely many
geometrically different soliton deformations of
a torus in $\R^4$ described by
the same DS equation, these different curves are just the values of the same
conservation law for different solutions.

\end{document}